\documentclass{ifacconf}

\usepackage{graphicx}      
\usepackage{natbib}        

\usepackage[utf8]{inputenc}
\usepackage[T1]{fontenc}
\usepackage[english]{babel}
\usepackage{amsmath}
\usepackage{amssymb}
\usepackage{amsfonts}
\usepackage{gensymb}
\usepackage{calc}
\usepackage{empheq}
\usepackage{stackengine}
\usepackage{float}
\usepackage{stmaryrd}
\usepackage{multirow}
\usepackage{wrapfig}
\usepackage{bm}
\usepackage{tikz}
\usepackage[outdir=./]{epstopdf}
\usepackage{subfigure}
\newtheorem{theorem}{\textbf{Theorem}}

\newtheorem{remark}{\textbf{Remark}}

\newtheorem{problem}{\textbf{Problem}}

\DeclareMathOperator{\tr}{tr}
\makeatletter
\renewcommand*\env@matrix[1][*\c@MaxMatrixCols c]{%
	\hskip -\arraycolsep
	\let\@ifnextchar\new@ifnextchar
	\array{#1}}
\makeatother
\bibpunct{[}{]}{,}{a}{}{;}

\begin{document}
\begin{frontmatter}

\title{Optimal exploration strategies for finite horizon
 regret minimization in some adaptive control problems} 

\thanks[footnoteinfo]{This work was supported by VINNOVA Competence Center AdBIOPRO, contract [2016-05181]  and by the Swedish Research Council through the research environment NewLEADS (New Directions in Learning Dynamical Systems), contract [2016-06079] and by Wallenberg AI, Autonomous Systems and Software Program (WASP), funded by Knut and Alice Wallenberg Foundation. {\textit{E-mail adresses}}: K\'evin Colin: kcolin@kth.se, Håkan Hjalmarsson: hjalmars@kth.se, Xavier Bombois: xavier.bombois@ec-lyon.fr.}

\author[DCS]{K\'evin Colin}, 
\author[DCS]{H\r{a}kan Hjalmarsson},
\author[Ampere,CNRS]{Xavier Bombois} 

\address[DCS]{Division of Decision and Control Systems, KTH Royal Institute of Technology, Sweden}
\address[Ampere]{Laboratoire Ampère, UMR CNRS 5005, Ecole Centrale de Lyon, Université de Lyon, France}
\address[CNRS]{Centre National de la Recherche Scientifique (CNRS), France}

\begin{abstract}                
:\ In this work, we consider the problem of regret minimization in adaptive minimum variance and linear quadratic control problems. Regret minimization has been extensively studied in the literature for both types of adaptive control problems. Most of these works give results of the optimal rate of the regret in the asymptotic regime. In the minimum variance case, the optimal asymptotic rate for the regret is $\log(T)$ which can be reached without any additional external excitation. On the contrary, for most adaptive linear quadratic problems, it is necessary to add an external excitation in order to get the optimal asymptotic rate of $\sqrt{T}$. In this paper, we will actually show from an a theoretical study, as well as, in simulations that when the control horizon is pre-specified a lower regret can be obtained with either no external excitation or a new exploration type termed immediate. 
\end{abstract}


\end{frontmatter}

\section{Introduction}

Minimum variance (MV) controllers~\citep{Astrom:70} and linear quadratic regulators (LQR)~\citep{anderson2007optimal} are two examples of control policies for linear time invariant (LTI) systems for which exact knowledge of the dynamics is necessary to guarantee optimal control performances.

However, perfect knowledge of the dynamics (transfer function or state-space representation for LTI systems) is never possible to get due to the presence of unmeasured disturbances. A remedy to this problem is to implement an adaptive control law where the controller is updated online in order to compensate for performance losses.

One approach to adaptive control is known as \textit{optimism in face of uncertainty} (OFU) or \textit{bet on the best} approach where early work can be found in~\citep{Lai&Robbins:85} and then later in, e.g.,~\citep{pmlr-v19-abbasi-yadkori11a,Cam:97,lale2020explore}. However, such methods often lead to non-convex optimization problems. More recently, Thompson sampling has received significant attention for its empirical performances with low computational requirement. Such a method was recently applied to the adaptive LQR problem in~\citep{Abe:17,Ouy:17}.

Another class of adaptive controller strategies is based on the certainty equivalence principle~\citep{aastrom2013adaptive,Lai:86,Rantzer:18,Shirani:20,pmlr-v119-simchowitz20a,Cassel:20,jedra2021minimal}. In this case, the dynamics are recursively estimated by using all the data available since the beginning of the experiment and the controller is updated online pretending these estimates are exact. 

In many adaptive control problems, it is vital to add an external excitation in order to guarantee sufficient richness of the data and/or an appropriate decrease in the uncertainties of the identified models. However, this external excitation disturbs both the system output and the control effort which subsequently decreases the control performances. In both reinforcement learning and adaptive controller communities, significant effort has been spent in developing a framework in order to find an optimal trade-off between the performance degradation due to the uncertainties (exploitation cost) and the performance degradation due to the external excitation (exploration cost). It is called regret minimization, where the regret is a function of both the exploration and exploitation costs and the external excitation is designed in such a way that it minimizes the regret over an infinite or finite time horizon.

Early work on regret minimization can be found in~\citep{Lai:86,Lai:87_self} for MV controllers applied to single inpout single output (SISO) ARX systems. The regret is at best growing as $\mathcal{O}(\text{log}(T))$ asymptotically. Such a rate can be obtained without external excitation. 

More recent works on regret minimization focus on the adaptive LQR problem. The study of regret minimization in LQR problems was re-encouraged by the work in~\citep{pmlr-v19-abbasi-yadkori11a} where an OFU algorithm was developed providing a rate of $\mathcal{O}(\sqrt{T})$ for the regret. The robust controller design algorithm in~\citep{Dea:18} gives a regret asymptotically growing as $\mathcal{O}(T^{2/3})$. This rate was also guaranteed with a Thompson sampling approach in~\citep{Abe:17}. This was later improved in~\citep{Ouy:17} providing a rate of $\mathcal{O}(\sqrt{T})$. Further works established the same rate of $\mathcal{O}(\sqrt{T})$ for LQR problems with the certainty equivalence principle~\citep{Faradonbeh:19,Mania2019CertaintyEI}. The rate of $\mathcal{O}(\sqrt{T})$ is actually the optimal one for LQR problem and this was proven in~\citep{pmlr-v119-simchowitz20a} when both state matrices are unknown. Such rate can be achieved by exciting the system with a white Gaussian noise excitation whose variance decays as $\mathcal{O}(1/\sqrt{T})$~\citep{Wan:21}. In~\citep{jedra2021minimal}, it is proven that regret can be upper-bounded as $\mathcal{O}(\sqrt{T})$ for any time-horizon when both state matrices are unknown. 

From the rich literature on regret minimization, it seems that this subject has been well explored for the MV and LQR cases. However, in this paper we argue that for finite horizon problems immediate exploration may be better. This strategy entails to explore the system dynamics as early as possible after which only exploitation takes place. As grounds for this we provide analysis of an approximative model of the regret, inspired from application oriented experiment design~\citep{bombois:2006,jansson:2004,hjalmarsson:2009}. Simulation results supporting the conclusions from this analysis are also provided. The used regret model decouples the exploration and exploitation costs with both terms depending on the Fisher information matrix. From this model, we show that the optimal exploration strategy minimizing the regret in finite time is an immediate exploration for many MV and LQR problems. Some simulations are presented in order to support this new result.

\textbf{Notation.} The set of real-valued matrices of dimension $n\times m$ will be denoted $\mathbf{R}^{n\times m}$. When $\mathbf{A}$ is positive definite (resp. positive semi-definite), we will write $\mathbf{A}\succ 0$ (resp. $\mathbf{A}\succeq 0$). The expectation operator will be denoted by $\mathbb{E}$. The identity matrix of dimension $n\times n$ will be denoted by $\mathbf{I}_n$ { and $\mathbf{A}^\top$ denotes the transpose of any matrix $\mathbf{A}$. The notation $X\sim N(a,b)$ refers to the random variable $X$ which is normally distributed with mean $a$ and variance $b$. The discrete-time forward operator is denoted by $z$ and $\omega$ denotes the angular rate.}

\section{System and control objectives}

\subsection{Considered system}

Consider a SISO discrete-time LTI system $\mathcal{S}$ given by 
\begin{equation}
    y(t) = G(z) u(t) + H(z) e(t)
\end{equation}where $y$ and $u$ are the output and input respectively (resp.), $e$ is a zero-mean white Gaussian noise of variance $\sigma_e^2$, $G(z)$ is a stable transfer function with an unit time-delay and $H(z)$ a stable, inversely stable transfer function which is also monic (i.e., $H_0(z=\infty) = 1$). 

\subsection{Control Scenario 1: minimum variance}\label{sec:control_scenario_1}

In a first scenario, we wish to put the system under minimum variance control which ideally guarantees that $y(t) = e(t) \ \forall t$. From~\citep{Astrom:70}, the ideal control policy is given by $u(t) = -K_{mv}(z) y(t)$ with 
\begin{equation}\label{eq:ideal_min_var_K}
K_{mv}(z) = -G^{-1}(z)(1-H(z))
\end{equation}for the case that the system $\mathcal{S}$ is minimum-phase

\subsection{Control Scenario 2: linear quadratic control}\label{sec:control_scenario_2}

For the second control scenario, we will consider a state representation of the system $\mathcal{S}$
{\begin{align}
    x(t) &= \mathbf{A}x(t-1) +  \mathbf{B}u(t-1) +  v_x(t)\\
    y(t) &= \mathbf{C}x(t)  +  v_y(t)\label{eq:measurement_eq}
\end{align}}where $x\in\mathbb{R}^n$ is the state-vector of dimension $n$, $\mathbf{A}\in\mathbb{R}^{n\times n}$, $ \mathbf{B}\in\mathbb{R}^{n\times 1}$ and $\mathbf{C}\in\mathbb{R}^{1\times n}$. The noises $v_x$ and $v_y$ are zero-mean white Gaussian noises depending on $e(t)$. \textit{We will assume that we measure the state vector $x$}. We can implement  a state feedback linear quadratic control $u(t) = - \mathbf{K}_{lqr} x(t)$ which minimizes the  infinite time horizon cost $J_\infty(u) = \underset{{T\rightarrow+\infty}}{\lim}  J_T(u,x)/T$ with\begin{equation}\label{eq:cost_function} J_T(u) = \sum_{t=1}^T x(t)^\top\mathbf{Q}x(t) + u(t)^\top\mathbf{R}u(t)
\end{equation}where $\mathbf{Q}\succ 0$ and $\mathbf{R}\succeq 0$. It is well known that $\mathbf{K}_{lqr} = (\mathbf{R}+\mathbf{B}^\top\mathbf{P}\mathbf{B}^\top)^{-1}\mathbf{B}^\top\mathbf{P}\mathbf{A}$ where $\mathbf{P}$ is the unique positive definite solution to the following discrete-time algebraic Riccati equation (DARE){\begin{align}\label{eq:ideal_dare}
   \mathbf{P} &= \mathbf{A}^\top \mathbf{P}\mathbf{A} + \mathbf{Q}- \mathbf{A}^\top\mathbf{P}\mathbf{B}(\mathbf{R}+\mathbf{B}^\top\mathbf{P}\mathbf{B})^{-1}\mathbf{B}^\top\mathbf{P}\mathbf{A}
\end{align}}

\section{Adaptive control}

In both aforementioned control scenarios, the controller depends on the true dynamics of $\mathcal{S}$ which are unfortunately unknown to us. As a remedy, we put $\mathcal{S}$ under a certainty equivalence adaptive feedback policy in both situations. 

\subsection{Minimum variance adaptive controller (MVAC)}\label{sec:MVAC}

For Control Scenario 1, the policy is given by
\begin{equation}
    u(t) = -\hat{K}_{mv}(t,z)y(t) + w(t) 
\end{equation}where $\hat{K}_{mv}(t,z)$ is the certainty equivalence controller transfer function and $w_{mv}$ an external excitation used to gather more information about $\mathcal{S}$. The adapted controller $\hat{K}_{mv}(t,z)$ is given by $ \hat{K}_{mv}(t,z) = -\hat{G}^{-1}(t,z)(1-\hat{H}(t,z))$ where $\hat{G}(t,z)$ and $\hat{H}(t,z)$ are respectively the identified transfer function of $G(z)$ and ${H}(z)$ obtained at time instant $t$ using, e.g., prediction error identification with input-output data $\{u(k),y(k)\}_{k=1}^t$ up to time instant $t$.

\subsection{Linear quadratic adaptive controller (LQAC)}\label{sec:LQAC}

In Control Scenario 2, the policy is given by
\begin{equation}
    u(t) = -\hat{\mathbf{K}}_{lqr}(t)x(t) + w(t) 
\end{equation}where $\hat{\mathbf{K}}_{lqr}(t) = (\mathbf{R}+\hat{\mathbf{B}}(t)^\top\hat{\mathbf{P}}(t)\hat{\mathbf{B}}(t)^\top)^{-1}\hat{\mathbf{B}}(t)^\top\hat{\mathbf{P}}(t)\hat{\mathbf{A}}(t)$ is the certainty equivalence controller, $\hat{\mathbf{A}}(t)$ and $\hat{\mathbf{B}}(t)$ are respectively the identified state matrices of $\mathbf{A}$ and $\mathbf{B}$ obtained at time instant $t$ using least-squares identification with input-state data\footnote{Recall that we assume that we measure the state vector $x$.} $\{u(k),x(k)\}_{k=1}^t$ up to time instant $t$. The matrix $\hat{\mathbf{P}}(t)$ is the positive definite of the DARE in~\eqref{eq:ideal_dare} for which $\mathbf{A}$ and $\mathbf{B}$ are respectively replaced by $\hat{\mathbf{A}}(t)$ and $\hat{\mathbf{B}}(t)$.

\section{Regret and previous results}

\subsection{Regret minimization}

In the recent adaptive control literature, much of the work focuses on designing the external excitation $w$ in order to reach an optimal trade-off between exploration (actions applied on the system for gathering information, in our case it is the signal $w$) and exploitation (actions applied on the system for control cost minimization, i.e., the controllers  $\hat{K}_{mv}(t,z)$ and $\hat{\mathbf{K}}_{lqr}(t)$). The literature often tackles the problem of trade-off by defining a function called cumulative regret $R(T)$. The optimal trade-off is then obtained by designing the sequence $\{w(t)\}_{t=1}^T$ such that $R(T)$ is minimized or such that the growth rate of the regret is minimized.

\subsection{Regret for MVAC}

Early work of regret minimization can be found in~\citep{Lai:87_self} where the MVAC of SISO ARX systems is treated. In this work, the regret is defined as 
{\begin{equation}\label{eq:regret_mv}
    R_{mv}(T) = \sum_{t=1}^T \mathbb{E}[(y(t) - e(t))^2]
\end{equation}}where $y$ is the output of $\mathcal{S}$ put under MVAC. We have the following result from the literature: 
\begin{theorem}[\citep{Lai:87_self}]\label{theo:MVAC}
    \textit{Consider the framework of the MVAC of Control Scenario 1 defined above. The optimal asymptotic rate for the regret $R_{mv}(T)$ is $\log(T)$. Such optimal rate can be reached without external excitation, i.e., $w = 0$.}
\end{theorem}We will refer to the type of exploration $w = 0$ as lazy.

\subsection{Regret for LQAC}

Most of the recent work has focused on the regret for LQAC. In~\citep{Wan:21}, the regret minimized is the following one
\begin{equation}\label{eq:regret_lqr}
    R_{lqr}(T) = \mathbb{E}[J_T(u)] - \mathbb{E}[J_T(u^*)] 
\end{equation}where $u^*$ is the input when the optimal linear quadratic controller $\mathbf{K}_{lqr}$ is applied to $\mathcal{S}$. We have the following result from the literature~\citep{Wan:21}:
\begin{theorem}\label{theo:LQAC}
    \textit{Consider the framework of the LQAC of Control Scenario 2. When both state matrices $\mathbf{A}$ and $\mathbf{B}$ are unknown, the optimal asymptotic rate for the regret $R_{lqr}(T)$ is $\sqrt{T}$. Such optimal rate can be reached with a white Gaussian noise external excitation of the form
    {\begin{equation}\label{eq:decaying_exploration}
        w(t)\sim N\left(0,\alpha/\sqrt{t} \right) \ \ \ \  \alpha \ge 0
    \end{equation}}We will refer to this type of exploration as $1/\sqrt{t}$-decaying.}  
\end{theorem}We will refer to the type of exploration~\eqref{eq:decaying_exploration} as $1/\sqrt{t}$-decaying.
Another way to achieve this optimal rate of Theorem~\ref{theo:MVAC} is to consider Thompson sampling~\citep{Ouy:17,Abe:17}. Since it is not based on the certainty equivalence principle, we will not consider Thompson sampling in this paper.

The results available in the literature give the optimal rate in the asymptotic regime and an upper bound of $\mathcal{O}(\sqrt{T})$ in~\citep{jedra2021minimal} when both state matrices $\mathbf{A}$ and $\mathbf{B}$ are unknown. In the following section we will examine the finite horizon case, i.e., when the cumulative regret over a fixed finite horizon $T$ is considered. We will assume $T$ is large so that asymptotic arguments can  still be used. 

\section{Abstract theoretical study}

\subsection{A novel regret model}

Denote by $\theta$ the vector of model parameters of dimension $n_\theta$ and $\theta_0$ the true parameter vector. Regret minimization can be seen as an experiment design problem since we look for the optimal sequence $\{w(t)\}_{t=1}^T$ such that a given criterion (here $R(T)$) is minimized. In the classical literature of application-oriented experiment design of LTI systems~\citep{bombois:2006,jansson:2004,hjalmarsson:2009}, we design an optimal finite-time identification experiment such that the expected value of the performance degradation $D_T(\hat{\theta}_T,\theta_0)$ due to the parameter uncertainties of $\hat{\theta}_T$ is minimized under some constraints on the input power. By using a Taylor expansion up to the second degree, we have $\mathbb{E}[D_T(\hat{\theta}_T,\theta_0)] \approx   \mathbb{E}[(\hat{\theta}_T-\theta_0)^\top \mathbf{W}(\theta_0)(\hat{\theta}_T-\theta_0)] = \tr(\mathbf{W}(\theta_0) \mathbb{E}[(\hat{\theta}_T-\theta_0)(\hat{\theta}_T-\theta_0)^\top])$ where $2\mathbf{W}(\theta_0)$ is the Hessian of $D_T(\hat{\theta}_T,\theta_0)$ evaluated at $\hat{\theta}_T = \theta_0$. Assuming an efficient estimator, we have $\mathbb{E}[(\hat{\theta}_T-\theta_0)(\hat{\theta}_T-\theta_0)^\top] = \mathcal{I}_{T}^{-1}$ where $\mathcal{I}_{t}$ is the Fisher information matrix. When the noise is Gaussian, it is given by
{\small\begin{equation}\label{eq:fisher}
   \mathcal{I}_{t} = \mathcal{I}_{t-1} + \mathcal{L}_t \ \ \text{with} \ \ \mathcal{L}_t = \dfrac{1}{\sigma_e^2}\mathbb{E}\left[\nabla_{\theta}\hat{\epsilon}(t,\theta)\left(\nabla_{\theta}\hat{\epsilon}(t,\theta)\right)^\top\right]|_{\theta = \theta_0}
\end{equation}}where $\nabla_\theta$ is the gradient operator w.r.t. $\theta$ and where $\epsilon(t,\theta)$ is the prediction error. We will write $\nabla_{\theta}\epsilon(t,\theta)$ as follows $\nabla_{\theta}\epsilon(t,\theta) = F_e(z,\theta)e(t) + F_w(z,\theta)w(t)$ with $F_e(z,\theta)$ and $F_w(z,\theta)$ depending on the model structure and the closed-loop transfer functions. In the regret minimization case, we have two types of performance degradation costs: the cumulative performance degradation cost due to the uncertainties and the cumulative performance  degradation cost due to the external excitation. The former can be modeled by the sum of all terms $\mathbb{E}[D_T(\hat{\theta}_t,\theta_0)] \approx \tr(\mathbf{W}(\theta_0) \mathcal{I}_{t}^{-1})$ from $t=1$ till $t=T$. We will model the latter with the sum of the power $\mathcal{P}_{w,t}$ of the signal $w$ at time $t$. Finally, by assuming that these costs are additive, we get\vspace{-0.2cm}
{\small\begin{equation}\label{eq:interm}
    R(T) \approx  \sum_{t=1}^T \left(\tr(\mathbf{W}(\theta_0)\mathcal{I}_{t}^{-1}) +  \mathcal{P}_{w,t}\right)\vspace{-0.2cm}
\end{equation}}With the independence assumption between $e$ and $w$, $\mathcal{L}_t$ in~\eqref{eq:fisher} can be split into two terms $\mathcal{L}_{t} = \mathcal{L}_{e,t} + \mathcal{L}_{w,t}$ where $\mathcal{L}_{e,t} = 1/\sigma_e^2\mathbb{E}[(F_e(z,\theta_0)   e(t))(F_e(t,\theta_0)  e(t))^\top ]$ and $
    \mathcal{L}_{w,t} = 1/\sigma_e^2\mathbb{E}[(F_w(z,\theta_0)  w(t))(F_w(z,\theta_0)   w(t))^\top]$ Assume that $w(t) = \sigma_w(t) w_0(t)$ with $w_0(t)$ a zero-mean filtered white noise with unit variance and $\sigma_w(t)$ deterministic and varying slower than the time constant of $F_w(z,\theta_0)$. Then,  $\mathcal{L}_{w,t} \approx \sigma_w^2(t)/\sigma_e^2 \mathbb{E}[ (F_w(z,\theta_0)  w_0(t))(F_w(z,\theta_0)   w_0(t))^\top]$. The power $\mathcal{P}_{w,t}$ of $w(t)$ is $\sigma_w^2(t)$ which leads to $\mathcal{L}_{w,t} = \mathcal{P}_{w,t}/\sigma_e^2 \mathbb{E}[(F_w(z,\theta_0)  w_0(t))(F_w(z,\theta_0)   w_0(t))^\top]$. By introducing  $\mathbf{Z} = n_\theta/\sigma_e^2 \mathbb{E}[(F_w(z,\theta_0)  w(t))(F_w(z,\theta_0)   w(t))^\top]$ and assuming it is invertible, we have $\mathbf{Z}^{-1}\mathcal{L}_{w,t} = \mathcal{P}_{w,t}/n_\theta \mathbf{I}_{n_\theta}$. Taking the trace of the latter, we get $\mathcal{P}_{w,t} = \text{tr}(\mathbf{Z}^{-1}\mathcal{L}_{w,t})$. Injecting that into~\eqref{eq:interm}, we get the following model
{\small\begin{align}
       R(T)  &  \approx  \sum_{t=1}^T\left( \tr(\mathbf{W}\mathcal{I}_t^{-1}) +  \text{tr}(\mathbf{Z}^{-1}\mathcal{L}_{w,t})\right)\label{eq:model_regret}\\
        \mathcal{I}_t & =
 \mathcal{I}_{t-1} + \mathcal{L}_{e,t} + \mathcal{L}_{w,t}\label{eq:model_fisher}
\end{align}}
\begin{remark}
    \textit{Even though the model in~\eqref{eq:model_regret}-\eqref{eq:model_fisher} comes from several approximations, a similar expression has actually been proven for the LQR problem in Chapter 6 of~\citep{ferizbegovic:2022}. Since the minimum variance problem is a particular case of the LQR problem, this decoupling also holds in that case.}
\end{remark}

The instantaneous information matrices $\mathcal{L}_{e,t}$ and $\mathcal{L}_{w,t}$ are structured and depend on the true parameter vector $\theta_0$~\citep{Ljung:1999}. This means that minimizing the regret, as given by~\eqref{eq:model_regret}, subject to the evolution~\eqref{eq:model_fisher} of the Fisher matrix, with respect to (w.r.t.) the external excitation $\{w(t)\}$ is a hard task. However, by relaxing the problem so that  $\mathcal{L}_{w,t} $ is simply a positive semi-definite matrix we obtain a lower bound on the achievable regret. In this case, we can see  $\mathcal{L}_{w,t} $ as the decision variables, i.e. they represent the actions taken for exploring the system dynamics. To simplify matters we will consider the case where  $\mathcal{L}_{e,t}=\mathcal{L}_e $ is time-invariant. Therefore, from now on, we will focus on studying the closed-form solutions of the following problem
\begin{problem}\label{prob:unstr}
   {\small \begin{align*}
        & \ \ \underset{\mathcal{L}_{w,t} \succeq 0}{\text{ min }}  \sum_{t=1}^T\left( \tr(\mathbf{W}\mathcal{I}_t^{-1}) +  \text{tr}(\mathbf{Z}^{-1}\mathcal{L}_{w,t})\right)  \text{ with }    \mathcal{I}_t = t\mathcal{L}_e + \sum_{k=1}^t \mathcal{L}_{w,k}
    \end{align*}}
\end{problem}

\subsection{Solution for Problem 1 and interpretation}

The following theorem gives the main result of the paper.
\begin{theorem}\label{theo:result_paper}
    \textit{Consider Problem 1. Decompose $\mathbf{Z}^{-1}$ 
as follows $\mathbf{Z}^{-1} = \mathbf{V}\mathbf{V}$ where $\mathbf{V}$ is the unique positive definite matrix square root. Define  $c_T$ as follows
   {\small \begin{equation}   \label{eq:ct}   c_T=\left(\lambda_{\max}\left(\mathbf{Z}{\mathcal{L}}_e^{-1}{\mathbf{W}}{\mathcal{L}}_e^{-1}\right) \sum_{t=1}^T\dfrac{1}{t^2}\right)^{-1}
    \end{equation}}where $\lambda_{\max}(S)$ denotes the maximal eigenvalue of any matrix $S$.
    There are two possible solutions for Problem~1:
    \begin{itemize}
        \item Solution 1: $\mathcal{L}_{w,t} = 0 \ \forall t \ge 1$.
        \item Solution 2: $\mathcal{L}_{w,1} \ne 0$ and $\mathcal{L}_{w,t} = 0 \ \forall t\ge 2 $.
    \end{itemize}The type of solution depends on the full-rankness of the matrix ${\mathcal{L}}_e$ and the validity of the following inequality}
     {\small\begin{equation}\label{eq:ineq}
      c_T \ge 1
    \end{equation}}\textit{We have the three cases
    \begin{itemize}
        \item Case 1: if ${\mathcal{L}}_e \succ 0$ and~\eqref{eq:ineq} is valid, then Solution~1 is optimal and the regret is given by
       {\small \begin{equation*}
            R(T) = \tr({\mathcal{L}}_e^{-1}{\mathbf{W}}) \sum_{t=1}^T \dfrac{1}{t} 
        \end{equation*}}
        \item Case 2: if ${\mathcal{L}}_e \succ 0$ but~\eqref{eq:ineq} is not valid, then Solution~2 is optimal and the regret is bounded as follows
        {\small\begin{equation*}
         (2-c_T)c_T\tr({\mathcal{L}}_e^{-1}{\mathbf{W}})  \sum_{t=1}^T \dfrac{1}{t}  \le R(T) \le \tr({\mathcal{L}}_e^{-1}{\mathbf{W}}) \sum_{t=1}^T \dfrac{1}{t} 
        \end{equation*}}
        \item Case 3: if $\mathcal{L}_e$ is singular, then Solution~2 is optimal and the  regret is bounded as follows
        {\small\begin{equation*}
         \delta \sqrt{T} \le R(T) \le \chi \sqrt{T} + \nu_T
        \end{equation*}}where $\delta$, $ \chi$ and $\nu_T = o(\sqrt{T})$ are positive scalars whose expressions are available in the appendix enclosed with this paper.
    \end{itemize}}
\end{theorem}
\begin{pf}
    See the appendix enclosed with this paper.
\end{pf} Solution 1 corresponds to $w(t) = 0 \ \forall t\ge1$ which is lazy exploration. Solution 2 means that we only excite at the first time instant (i.e., $w(1) \ne 0$ and $w(t) = 0 \ \forall t\ge 2$) and we will call this type of excitation {\it immediate}. Notice that $c_T$ in~\eqref{eq:ct} will grow as the time horizon $T$ shrinks or as the quantity ${\mathcal{L}}_e^{-1}{\mathbf{W}}{\mathcal{L}}_e^{-1}$ decreases. This means that the results of Theorem~3 are intuitively appealing: it does not pay off to explore if the time horizon is too short and/or if the benefits from the noise excitation is sufficiently large in comparison to the regret incurred by the model error.

Several interesting connections can be made between the results in Theorem~\ref{theo:result_paper} and existing ones in the literature. First of all, we get the optimal rate of $\log(T)$ when $T\rightarrow +\infty$ with lazy exploration (Case 1 of Theorem~\ref{theo:result_paper} with Solution 1) as was the case for MVAC in Theorem~\ref{theo:MVAC}. Secondly, when $T\rightarrow +\infty$, we get the rate of $\sqrt{T}$ in Case 3 which is shown to be optimal for LQAC in Theorem~\ref{theo:LQAC}. This suggests that, although approximate in nature, Problem 1 may have bearings on actual adaptive control problems. To further examine this connection notice that, in the MVAC case, Theorem~\ref{theo:result_paper} suggests that there may be cases for which an immediate exploration would be better than a lazy exploration (depending if~\eqref{eq:ineq} is valid or not). Also, the optimal exploration we get in Case 3 for a rate of $\sqrt{T}$ is immediate and not an excitation with a variance decaying as $1/\sqrt{t}$ as proposed in the literature for the LQAC case (see Theorem~\ref{theo:LQAC}). In the next section, we examine in a numerical example if these  results for Problem~1 carry over to adaptive control problems.

\section{Numerical example}

\subsection{System}
We will consider the following first order ARX system
\begin{equation}\label{eq:num_ex_tf}
   ({1+a_0z^{-1}}) y(t) = b_0z^{-1} u(t) +  e(t)
\end{equation} We choose $a_0 = -0.45$ and $b_0 = 0.67$. We would like to investigate if there are cases for which immediate exploration is better than lazy exploration for MVAC and $1/\sqrt{t}$-decaying exploration for LQAC. For the model structure, we consider the same order as the true system. In both cases, the identification is a linear least squares problem so we can use recursive linear least-squares identification (see, e.g., Chapter 11 of~\citep{Ljung:1999}). 

\subsection{Initial phases and immediate exploration}

For lazy exploration, as proposed in~\citep{Lai:87_self}, we start from an initial covariance matrix $\mathbf{P}_{init}$ and initial parameters $a_{init}$ and $b_{init}$ for $a_0$ and $b_0$ respectively. The initial controller is built from $a_{init}$ and $b_{init}$. We choose $\mathbf{P}_{init} = 10^3\mathbf{I}_2$, $a_{init} = -0.3$ and $b_{init} = 0.8$. For $1/\sqrt{t}$-decaying exploration, $w$ is equal to~\eqref{eq:decaying_exploration} and we consider an open-loop initial phase of duration $N_i = 3$. We then compute the model estimate at $t=N_i$ and for $t\ge N_i$ we use the adaptive controller with certainty equivalence principle. For immediate exploration, the theory suggests that we only use the first time instant in order to compute an initial model and covariance matrix. However, the covariance matrix will never be positive definite in that case which is an issue for the recursive least-squares identification. Therefore, for immediate exploration, the initial phase is identical as for $1/\sqrt{t}$-decaying exploration except that $w$ is equal to a constant $\beta$ at $t=1$ then it is taken equal to $0$ afterwards. The certainty equivalence controller is designed  based on the model at $t=N_i$ and then \textit{it stays constant for the rest of the experiment}.

\subsection{MVAC: lazy or immediate?}\label{sec:example_MVAC}

Consider Control Scenario 1 in Section~\ref{sec:control_scenario_1} with the system described by~\eqref{eq:num_ex_tf} and controlled by the MVAC given in Section~\ref{sec:MVAC}. We know that $\log(T)$ is the optimal asymptotic rate for the regret $R_{mv}(T)$ in~\eqref{eq:regret_mv}. In our theory, this rate is achieved in Cases 1 and 2 of Theorem~\ref{theo:result_paper} where the only difference is the validity of the inequality~\eqref{eq:ineq}. In this case, the performance degradation is $\sum_{t=1}^T\mathbb{E}[(y(t) - e(t))^2]$ which increases with the noise variance $\sigma_e^2$.  Therefore, the Hessian $2\mathbf{W}$ increases when $\sigma_e^2$ increases. Hence, from our theory, by changing $\sigma_e^2$, we can change the validity of the inequality~\eqref{eq:ineq} which in turn changes the type of optimal exploration of Problem~\ref{prob:unstr} (immediate or lazy exploration). 

We will now examine if this holds also for MVAC in a simulation study by considering a gridding of $\sigma_e^2$  with $100$ log-regularly spaced values between $10^{-5}$ and $10^0$. We choose $T = 10^5$ for the time horizon. We  store $1000$ different realizations of a zero-mean white noise $\bar{e}(t)$ with \textit{unit variance}. Denote by $\bar{e}_i(t)$ the $i$-th stored noise realization of $\bar{e}(t)$ for the $i$-th simulation. For each value of $\sigma_e^2$, we simulate the MVAC with a lazy exploration 1000 times, one time per stored noise realization $\bar{e}_i(t)$, by setting the noise $e(t)$ equal to ${e}(t) = \sigma_e \bar{e}_i(t)$. For the immediate exploration, we do similarly but we need to determine the optimal immediate exploration gain $\beta$ for each value of $\sigma_e^2$. For this purpose, we do a gridding of $\beta$ with $100$ log-regularly spaced values between $10^{-5}$ and $10^0$. For each value of $\sigma_e^2$, we repeat the following procedure: we select one value of $\beta$ in the gridding and we simulate the MVAC $1000$ times, one per stored noise realization of $\bar{e}_i(t)$, by setting the noise $e(t)$ equal to ${e}(t) = \sigma_e \bar{e}_i(t)$  for the $i$-th simulation. Then, we approximate the regret with the average over the 1000 obtained values. We repeat the latter for every value of $\beta$ and after that, we select the $\beta$ giving the minimal regret among the 100 different cases and use the corresponding regret for the immediate exploration case. We switch to the next value of $\sigma_e^2$ in the considered gridding and repeat this procedure.

 In Figure~1.(a), we depict the cumulative regret $R_{mv}(T)$ at time instant $T = 10^5$ obtained with lazy exploration (blue dashed line) and the optimal immediate exploration (orange solid line)  for each value of $\sigma_e^2$ considered in the gridding. We clearly see two different cases for which a switch is observed at around $\sigma_e^* = \sqrt{2 \times 10^{-3}}$. For $\sigma_e < \sigma_e^*$ (resp. $\sigma_e > \sigma_e^*$), lazy exploration (resp. immediate exploration) is optimal. This observation supports the results from our abstract theoretical study, provided in Theorem~\ref{theo:result_paper}. Note that for $\sigma_e > \sigma_e^*$, the results of lazy exploration are not smooth anymore despite the 1000 Monte Carlo simulations considered in the study. In Figure~1.(b), we depict 6 Monte Carlo simulations of the evolution of $\sum_{t=1}^T(y(t) - e(t))^2$ with respect to $\sigma_e^2$. We observe that they are smooth but some peaks can happen which explains the behavior of the regret for $\sigma_e > \sigma_e^*$. We would need many more Monte Carlo simulations to get a smoother curve which takes a very long time. This also shows that it may be risky to use lazy exploration for one experiment. Further study of this phenomenon seems warranted.

\begin{figure*}
    \centering
    \subfigure[]{\includegraphics[width=0.326\textwidth]{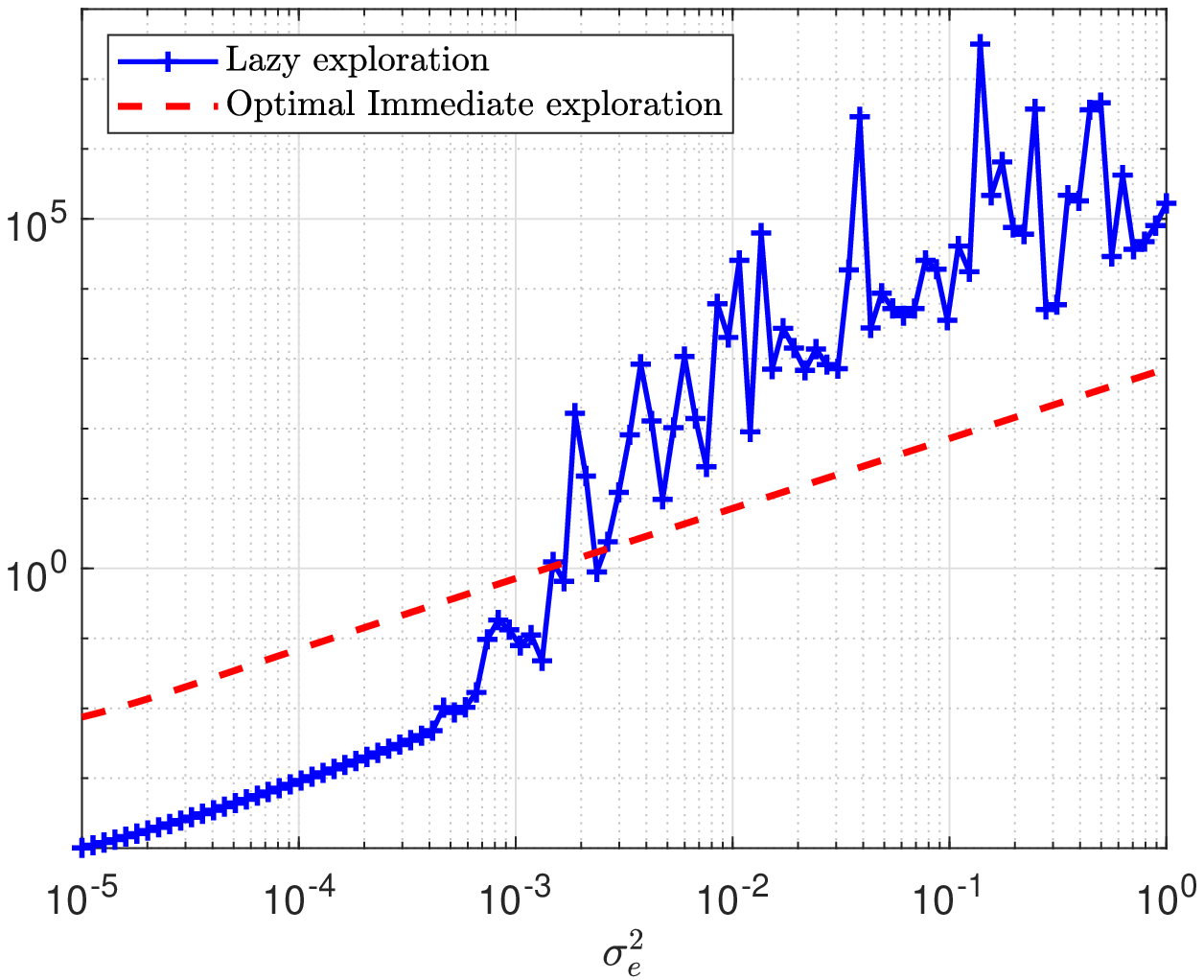}} 
    \subfigure[]{\includegraphics[width=0.326\textwidth]{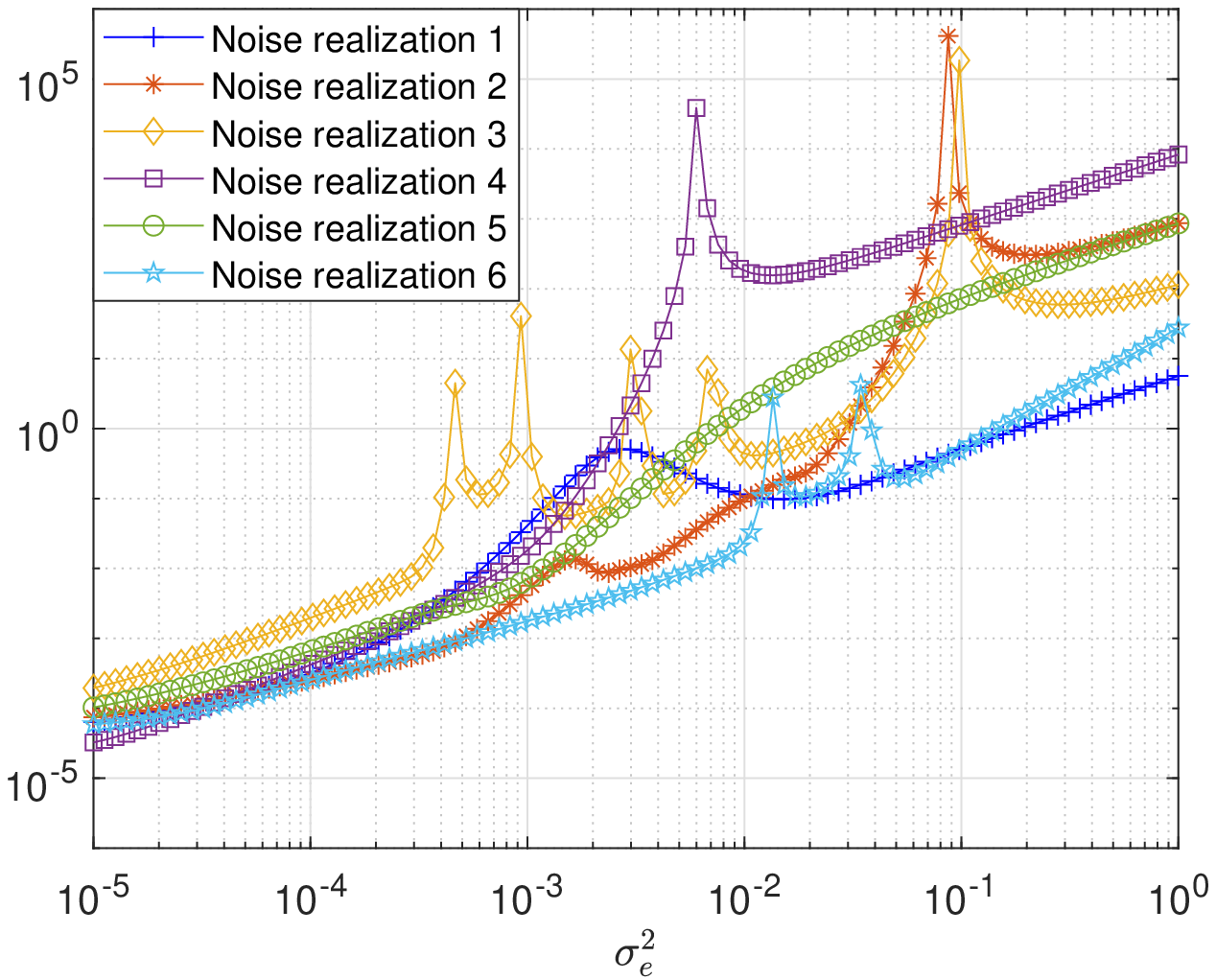}} 
    \subfigure[]{\includegraphics[width=0.326\textwidth]{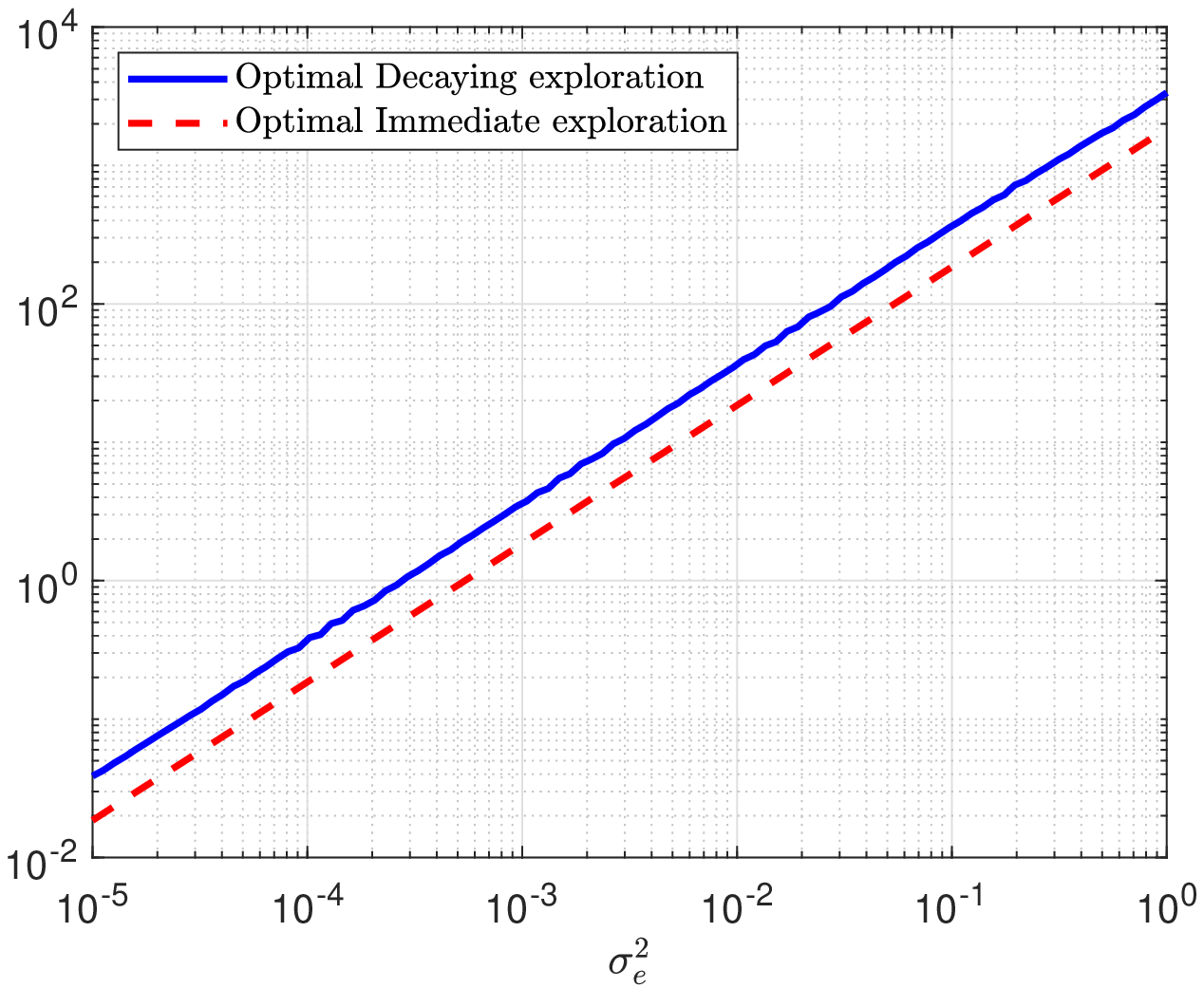}}
    \vspace{-0.35cm}\caption{ (a) Regret $R_{mv}(T)$ with MVAC at  $T = 10^5$ obtained with lazy exploration (blue line) and optimal immediate exploration (orange line) for different $\sigma_e^2$. (b) Value of $\sum_{t=1}^T(y(t) - e(t))^2$ for several noise realizations for MVAC  for different $\sigma_e^2$.  (c) Regret $R_{mv}(T)$ with LQAC at $T = 10^5$ obtained with optimal $1/\sqrt{t}$-decaying exploration (blue  line) and optimal immediate exploration (orange line)  for different $\sigma_e^2$.}
    \label{fig:all_figs}
\end{figure*}

\subsection{LQAC: $1/\sqrt{t}$-decaying or immediate?}

Consider Control Scenario 2 in Section~\ref{sec:control_scenario_2} with the system described by~\eqref{eq:num_ex_tf} and controlled by the LQAC given in Section~\ref{sec:LQAC}. We know that $\sqrt{T}$ is the optimal asymptotic rate for the regret $R_{lqr}(T)$ in~\eqref{eq:regret_lqr} and that it can be reached with a white Gaussian noise excitation with a variance decaying as $1/\sqrt{t}$ (see Theorem~\ref{theo:LQAC}). Theorem~\ref{theo:result_paper} implies that we may get lower regret with an immediate exploration. We do the same gridding for $\sigma_e^2$. We store 100 new noise realizations of the zero-mean white noise $\bar{e}(t)$ with \textit{unit variance} and denote by $\bar{e}_i(t)$ the $i$-th realization $(i=1,\cdots,100)$. For the immediate exploration case, we do the simulations as explained in Section~\ref{sec:example_MVAC} with the same gridding for $\beta$. For the white noise $1/\sqrt{t}$-decaying exploration, we store $100$ realizations of a zero-mean white noise $\bar{w}(t)$ with \textit{unit variance} and denote by $\bar{w}_i(t)$ the $i$-th realization $(i=1,\cdots,100)$. For each value of $\sigma_e^2$, we need to determine the optimal constant $\alpha$. As for $\beta$, we do a gridding of $\alpha$ with $100$ log-regularly spaced values between $10^{-5}$ and $10^0$. Then, we do the following procedure for each value of  $\sigma_e^2$: we select one value of $\alpha$ in the gridding and we simulate the LQAC $100$ times, by setting ${e}(t) = \sigma_e \bar{e}_i(t)$ and ${w}(t) = (\alpha/\sqrt{t})^{1/2} \bar{w}_i(t)$ for the $i$-th realization. Then, we approximate the regret with the average over the 100 obtained values. We repeat the latter for every value of $\alpha$ and after that, we select the minimal regret among the 100 computed values. We switch to the next value of $\sigma_e^2$ in the considered gridding and repeat this procedure. In Figure~1.(c), we depict the regret at time instant $T = 10^5$ obtained with the optimal $1/\sqrt{t}$-decaying exploration and the optimal immediate exploration for each value of $\sigma_e^2$ considered in the gridding. We observe that immediate exploration always performs  better than $1/\sqrt{t}$-decaying exploration which is consistent with the conclusions drawn in Theorem~\ref{theo:result_paper} for Problem~\ref{prob:unstr}.

\section{Discussion and conclusions}

We have observed new results on regret minimization: for a finite horizon $T$ we have shown in simulations that immediate excitation reduces the cumulative regret over the previously considered optimal policies of lazy exploration (MVAC) and $1/\sqrt{t}$-decaying exploration (LQAC). These results were predicted from the behaviour of the approximate regret minimization problem stated as Problem~\ref{prob:unstr}. as given in Theorem~\ref{theo:result_paper}. These observations tie in with the simulation results in~\citep{forgione:2015} on experiment design for adaptive $\mathcal{H}_2$ controllers which also show that the exploration effort should be distributed at the beginning of the experiment. The difference with our framework is that the controller is kept constant during a sufficient number of time instants before being updated so that the stationary assumption holds.

Our simulation results suggest that, in the case of an a priori known horizon $T$ for both LQAC and MVAC, using adaptive controllers may be useless since a constant controller obtained with immediate exploration reduces better the regret. However, the solution to the immediate exploration problem requires knowledge of  model parameters and noise variances. Moreover, the magnitude of the pulse of immediate exploration may be too large w.r.t. the system limitations. So in practice $1/\sqrt{t}$-decaying or lazy exploration may still be the preferred choice. Nevertheless, the observations made in this paper may serve as basis for reducing the regret in data driven control problems. In future works we will study ramifications of Problem 1.

\bibliography{ifacconf}             
                                                   







\onecolumn

\textbf{Appendix of the paper}

\section{Change of notations}

When $\mathbf{A}$ is positive definite (resp. positive semi-definite), we will write $\mathbf{A}\succ 0$ (resp. $\mathbf{A}\succeq 0$). The identity matrix will be denoted by $I$ { and $\mathbf{A}^\top$ denotes the transpose of any matrix $\mathbf{A}$. The trace operator is denoted by $\text{tr}\{ \}$.
We consider the following change in notations with respect to the paper: $R(T)\rightarrow R_T$, $\mathcal{I}_{t} \rightarrow \bar{X}_t$, $\mathcal{L}_{w,t} \rightarrow X_t$, $\mathcal{L}_{e} \rightarrow S$ and  $\mathbf{W} \rightarrow W$.

\section{Change of variables}\label{sec:change_of_variable}

The problem considered in the paper is the following one
 { \begin{align} \label{eq:prob_origi}
        & \ \ \underset{ X_t \succeq 0}{\text{ min }}  \underbrace{\sum_{t=1}^T\left( \text{tr}\left\{W\bar{X}_t^{-1}\right\} +  \text{tr}\left\{\mathbf{Z}^{-1}X_t\right\}\right)}_{R_T}  \text{ with }    \bar{X}_t = tS + \sum_{k=1}^tX_k
    \end{align}}Let us decompose $\mathbf{Z}^{-1}$ 
as follows $\mathbf{Z}^{-1} = \mathbf{V}\mathbf{V}$ where $\mathbf{V}$ is the unique positive definite matrix square root. Then, we can simplify~\eqref{eq:prob_origi} by introducing $X_{\mathbf{V},t} = \mathbf{V} X_t\mathbf{V}$ and $\bar{X}_{\mathbf{V},t} = \mathbf{V} \bar{X}_t\mathbf{V}$ so that
\begin{align}
\label{Itgenb2}
    \bar{X}_{\mathbf{V},t}&= \mathbf{V}\bar{X}_{0}\mathbf{V}+\sum_{k=1}^t X_{\mathbf{V},k}+t\mathbf{V}S\mathbf{V}=\bar{X}_{\mathbf{V},0}+\sum_{k=1}^t X_{\mathbf{V},k}+tS_{\mathbf{V}}\\
    \label{RTgenb}
    R_T&=\sum_{t=1}^T\text{tr}\left\{W(\mathbf{V}^{-1}\bar{X}_{\mathbf{V},t}\mathbf{V}^{-1})^{-1}\right\}+\sum_{t=1}^T\text{tr}\left\{\mathbf{Z}^{-1}\mathbf{V}^{-1}X_{\mathbf{V},t}\mathbf{V}^{-1}\right\}\\
    &=\sum_{t=1}^T\text{tr}\left\{W_\mathbf{V}\bar{X}_{\mathbf{V},t}^{-1}\right\}+\sum_{t=1}^T\text{tr}\left\{X_{\mathbf{V},t}\right\}
\end{align}
where $S_{\mathbf{V}}:=\mathbf{V}S\mathbf{V}$ and $W_{\mathbf{V}}:=\mathbf{V}W\mathbf{V}$. Thus, in the sequel, we will study the original problem~\eqref{eq:prob_origi} with $\mathbf{Z}^{-1} = I$. The solution of the original problem for any $\mathbf{Z}^{-1} = \mathbf{V}\mathbf{V}$ is obtained by doing the following change of variables: $\bar{X}_t \rightarrow \mathbf{V}\bar{X}_t\mathbf{V}$, $X_t \rightarrow \mathbf{V}X_t\mathbf{V}$, $S \rightarrow \mathbf{V}S\mathbf{V}$ and  $W \rightarrow \mathbf{V}W\mathbf{V}$.

\section{Theorem}\label{sec:theo}

\begin{theorem}
Let $X_t$ $(t=1,\ldots,T)$ be square matrices of the same dimensions and let $W\succ 0$, with eigen-decomposition $W=\sum_{k=1}^{n}\lambda_k e_ke_k^\top$,  and $S\succeq 0$ have the same dimensions as $X_t$. Let 
\begin{align}
    J_T(\{X_t\}_{t=1}^T):=\sum_{t=1}^T \tr\left\{W\bar{X}_t^{-1}\right\}+ \sum_{t=1}^T\tr\left\{X_t\right\}
\end{align}
where $ \bar{X}_t=\sum_{k=1}^TX_k+tS$.  

The problem
\begin{align}
\label{p5}
\begin{split}
    \inf_{X_t,t=1,\ldots,T} &  J_T(\{X_t\}_{t=1}^T)\\
    \text{s.t.} \; X_1+S&\succ0 \\
    X_t&\succeq 0,\; t=1,\ldots,T
    \end{split}
\end{align}
is convex with solution $\{X_t^*\}_{t=1}^T$ satisfying:
\begin{itemize}
    \item[i.] $X_t^*=0$, $t=2,\ldots,T$.
    \item[ii.] $X_1^*=0$ if and only if $S\succ0$ and 
    \begin{align}
    \label{lazycond}
     c_T:=\frac{1}{\beta_T\; \lambda_{\max}\left\{S^{-1}WS^{-1}\right\}}\geq 1  
    \end{align} 
    where $\beta_T:=\sum_{t=1}^T\frac{1}{t^2}$ and $\lambda_{\text{max}}\left\{X\right\}$ refers to the maximal eigenvalue of any matrix $X$, in which case
    \begin{align}
    \label{Joptii}
        J_T(\{X_t^*\}_{t=1}^T)=\alpha_T \;J_T^*
    \end{align}
    where $\alpha_T:=\sum_{t=1}^T\frac{1}{t}$ and $J_T^*:= \tr\left\{WS^{-1}\right\}$.

    Provided 
    \begin{align}
       I- \frac{\pi^2}{6}S^{-1}WS^{-1}\succ 0
    \end{align}
    it holds that
    \begin{align}
        \lim_{T\rightarrow\infty}\frac{ J_T(\{X_t^*\}_{t=1}^T)}{\log T}=J_T^*
    \end{align}
    \item[iii.] When $S\succ0$ but \eqref{lazycond} does not hold
    \begin{align}
        (2-c_T)c_T\leq \frac{ J_T(\{X_t^*\}_{t=1}^T)}{\alpha_T J_T^*}<1
    \end{align}
    \item[iv.] When $S$ is singular,
    \begin{align}
        &2\sqrt{T}\sum_{k=1}^mn_k^\top\tilde{W}n_k\leq  J_T(\{X_t^*\}_{t=1}^T)\\
            \label{SsingularJbounds}
        &\leq \sqrt{T}\left(\sum_{k=1}^m\frac{n_k^\top Wn_k}{n_k^\top\tilde{W}n_k}+\sum_{k=1}^mn_k^\top\tilde{W}n_k\right)+\gamma_T\sum_k \frac{e_k^\top We_k}{\lambda_k}
    \end{align}
    where $\gamma_T = \sum_{t=1}^T 1/t$,  $\{n_k\}_{k=1}^m$ is any orthonormal basis to the null space of $S$, where $\tilde{W}$ is the non-negative matrix square root of $W$ and where
    \begin{align}
        S=\sum_k \lambda_ke_ke_k^\top
    \end{align}
    is the eigen-decomposition of $S$.
    
    In particular when $S=0$,
    \begin{align}
    \label{sol5}
     X_1^*&=\sqrt{T} \tilde{W}\\
     \label{Jmin5}
      J_T(\{X_t^*\}_{t=1}^T)&= 2\sqrt{T}J_T^*
\end{align}
\end{itemize}
\end{theorem}

\begin{remark}
    Before dealing with the proof, we want to comment on the inequality~\eqref{lazycond} which is similar to the one of Theorem 3 of the paper but with $\mathbf{Z}^{-1} = I$. To consider any $\mathbf{Z}^{-1} = \mathbf{V}\mathbf{V}$, we first do the change of variables $S \rightarrow \mathbf{V}S\mathbf{V}$ and  $W  \rightarrow \mathbf{V}W\mathbf{V}$ (see Section~\ref{sec:change_of_variable}). Therefore, $\lambda_{\max}\left\{S^{-1}WS^{-1}\right\}$ is changed to
    \begin{align}
         \lambda_{\max}\left\{\mathbf{V}^{-1}S^{-1}\mathbf{V}^{-1}\mathbf{V}W\mathbf{V}\mathbf{V}^{-1}S^{-1}\mathbf{V}^{-1}\right\} &= \lambda_{\max}\left\{\mathbf{V}^{-1}S^{-1}WS^{-1}\mathbf{V}^{-1}\right\}
    \end{align}Now, by observing that $\mathbf{V}^{-1}S^{-1}WS^{-1}\mathbf{V}^{-1}$ and $\mathbf{V}^{-1}\mathbf{V}^{-1}S^{-1}WS^{-1}$ are similar matrices (they have the same eigenvalues), we have $\lambda_{\max}\left\{\mathbf{V}^{-1}S^{-1}WS^{-1}\mathbf{V}^{-1}\right\} = \lambda_{\max}\left\{\mathbf{V}^{-1}\mathbf{V}^{-1}S^{-1}WS^{-1}\right\}$. By recalling that $\mathbf{Z}^{-1} =\mathbf{V}\mathbf{V} $, we obtain $\lambda_{\max}\left\{\mathbf{Z}\mathbf{V}^{-1}S^{-1}WS^{-1}\right\}$ which is the expression in Theorem 3 of the paper.
\end{remark}

\section{Proof of Theorem~1 in Section~\ref{sec:theo}}
Introduce symmetric matrices $Z_t$, $t=1,\ldots,T$. Then the problem \eqref{p5} can be written as
\begin{align}
\inf_{Z_t,X_t,t=1,\ldots,T} & \sum_{t=1}^{T} \tr\left\{Z_t\right\}+\sum_{t=1}^T \tr\left\{X_t\right\}\\
\text{s.t.}\quad  \tilde{W}\bar{X}_t^{-1}\tilde{W}&<Z_t,\; t=1,\ldots,T\\
X_1+S&\succ0 \\
    X_t&\succeq 0,\; t=1,\ldots,T
\end{align}
Since $X_1+S\succ0$, $\bar{X}_t\succ0$ for $t=1,\ldots,T$ and therefore Schur complement (see, e.g., Appendix A.5.5 in \cite{Boyd&Vandenberghe:03}) gives that problem \eqref{p5} is equivalent to 
\begin{align}
\label{p6}
\inf_{Z_t,X_t,t=1,\ldots,T}  \sum_{t=1}^{T} \tr\left\{Z_t\right\}+\sum_{t=1}^T \tr\left\{X_t\right\}&\\
\text{s.t.}\quad  \begin{bmatrix}
Z_t &\tilde{W}\\\tilde{W} & \sum_{k=1}^tX_k+tS
\end{bmatrix}&\succeq 0,\; t=1,\ldots,T\\
X_1+S&\succ0\\
X_t&\succeq 0,\; t=1,\ldots,T
\end{align}
which is a semi-definite program (SDP) and therefore convex. 

We start with relaxing the strict constraint $X_1+S\succ0$ to a non-strict inequality, in which case we can remove the constraint since $S\succeq 0$ by assumption and $X_1\succeq0$ is included in the constraints in \eqref{p6}. 
With 
\begin{align}
    Q_t:=\begin{bmatrix}
    Q_{t,11} & Q_{t,12}\\
    Q_{t,12}^\top & Q_{t,22} 
    \end{bmatrix}\succeq 0,\; t=1,\ldots, T
\end{align} 
with each block $Q_{t,ij}$ having the same dimension as $X_t$, and 
$Y_t\succeq 0$, for $t=1,\ldots,T$, all having the same dimensions as $X_t$, and $Q=\{Q_{t}\}_{t=1}^{T}$ and  $Y=\{Y_{t}\}_{t=1}^{T}$,
the Lagrange dual function of this problem is given by 
\begin{align}
    {\mathcal L}(Q,Y)=&\inf_{Z_t,X_t,t=1\ldots,T} \sum_{t=1}^{T} \tr\left\{Z_t\right\}+\sum_{t=1}^T \tr\left\{X_t\right\}\\
    &-\sum_{t=1}^T\tr\left\{\begin{bmatrix}
    Q_{t,11} & Q_{t,12}\\
    Q_{t,12}^\top & Q_{t,22} 
    \end{bmatrix} \begin{bmatrix}
Z_t &\tilde{W}\\\tilde{W} & \sum_{k=1}^tX_k+tS
\end{bmatrix}\right\}-\sum_{t=1}^{T}\tr\left\{Y_tX_t\right\}\\
=&\inf_{Z_t,X_t,t=1\ldots,T} \sum_{t=1}^{T} \tr\left\{Z_t\right\}+\sum_{t=1}^T \tr\left\{X_t\right\}\\
&-\sum_{t=1}^T  \tr\left\{Q_{t,11}Z_t+2Q_{t,12}\tilde{W}+Q_{t,22}\left(\sum_{k=1}^tX_k+tS\right)\right\}-\sum_{t=1}^{T}\tr\left\{Y_tX_t\right\}\\
=&\inf_{Z_t,X_t,t=1\ldots,T} \sum_{t=1}^{T} \tr\left\{(I-Q_{t,11})Z_t\right\}+\sum_{t=1}^T \tr\left\{\left(I-\sum_{k=t}^TQ_{k,22}-Y_t\right)X_t\right\}\\
&-2\sum_{t=1}^T  \tr\left\{Q_{t,12}\tilde{W}\right\}-\sum_{t=1}^Tt \tr\left\{Q_{t,22}S\right\}
\end{align}
For this function to be finite valued, and not $-\infty$,
\begin{align}
\label{Qcond1}
Q_{t,11}&=I\\
\label{Qcond2}
    I-\sum_{k=t}^TQ_{k,22}&=Y_t\succeq 0,\quad t=1,\ldots,T
\end{align}
are required, in which case
\begin{align}
    {\mathcal L}(Q,Y)=&-2\sum_{t=1}^T  \tr\left\{Q_{t,12}\tilde{W}\right\}-\sum_{t=1}^Tt\tr\left\{Q_{t,22}S\right\}\\
    =&-2\sum_{t=1}^T  \tr\left\{Q_{t,12}\tilde{W}\right\}-\sum_{t=1}^Tt\tr\left\{Q_{t,22}S\right\}
\end{align}
Now according to duality theory (see, e.g., Section 5.1.3 in \citep{Boyd&Vandenberghe:03}), ${\mathcal L}(Q,Y)$ is a lower bound to objective function for the solution of the primal problem for any feasible $Q$, and $Y$. In fact as the original problem is strictly feasible, Slater's condition gives that maximum of ${\mathcal L}(Q,Y)$ equals $J_T(\{X_t^*\}_{t=1}^T)$. In view of \eqref{Qcond1}, Schur complement gives that
\begin{align}
    Q_t\succeq 0\quad \Leftrightarrow \quad Q_{t,22}\succeq Q_{t,12}^\top Q_{t,12}
\end{align}
Further, we note from that $Q_{t,22}\succeq 0$, that it is sufficient to require \eqref{Qcond2} for $t=1$.
We will thus study the problem 
\begin{align}
    \inf_{Q,Y} & -{\mathcal L}(Q,Y)\\
    \text{s.t } \quad Q_{t,22}&\succeq Q_{t,12}^\top Q_{t,12}\\
    I &\succeq \sum_{t=1}^TQ_{t,22}
\end{align}
In view of that the objective function is monotone increasing in $Q_{t,22}$ it follows that the optimum will satisfy $Q_{t,22}= Q_{t,12}^\top Q_{t,12}$. We can thus re-write the problem as
\begin{align}
\label{d5}
    \inf_{Q,Y}\quad & 2\sum_{t=1}^T  \tr\left\{Q_{t,12}\tilde{W}\right\}+\sum_{t=1}^Tt\tr\left\{Q_{t,12}SQ_{t,12}^\top\right\}\\
    \text{s.t }& I\succeq \sum_{k=1}^TQ_{t,12}^\top Q_{t,12}
\end{align}
This is a quadratically constrained quadratic program. We now consider the case where $S\succ0$ where we can express the (negative) dual function \eqref{d5} as
\begin{align}
    \sum_{t=1}^T  \tr\left\{\left(Q_{t,12}+\tilde{W}(tS)^{-1}\right)(tS)\left(Q_{t,12}+\tilde{W}(tS)^{-1}\right)\right\}-\tr\left\{\tilde{W}(tS)^{-1}\tilde{W}\right\}
\end{align}
whose unconstrained minimizer is given by $Q_{t,12}=-\frac{1}{t}\tilde{W}S^{-1}$, giving the objective function
\begin{align}
    -\sum_{t=1}^T\frac{1}{t} \tr\left\{\tilde{W}S^{-1}\tilde{W}\right\}=-\alpha_TJ_T^*
\end{align}
This is thus the optimum provided that the constraint in \eqref{d5} is met, i.e.
\begin{align}
   I\succ \sum_{t=}^T\frac{1}{t^2}S^{-1}WS^{-1}
\end{align}
which is equivalent to the constraint in \eqref{lazycond}. Summarizing, the maximum of the dual ${\mathcal L}(Q,Y)$ in the feasible set is given by $\alpha_T J_T^*$ when \eqref{lazycond} holds. Now, $X_t=0$, $t=1,\ldots,T$ is a feasible point for \eqref{p5} which has objective function $J_T(\{X_t\}_{t=1}^T)=\alpha_T J_T^*$, i.e. equal to the maximum of the dual. But by duality theory this must then be the solution to \eqref{p5}. This means that the if part in ii. of the theorem has been proven as well as \eqref{Joptii}. The final part of ii. follows by noticing that
\begin{align}
    \lim_{T\rightarrow\infty}\beta_T=\frac{\pi^2}{6}
\end{align}
For the cases when $S$ is singular or $S\succ0$ but \eqref{lazycond} does not hold, we proceed by computing the Lagrange dual function of \eqref{d5}, with $\Lambda\succeq 0$ as dual variable,
\begin{align}
   {\mathcal D}(\Lambda):=& \inf_{Q_{t,12},t=1,\ldots,T}\; 2\sum_{t=1}^T  \tr\left\{Q_{t,12}\tilde{W}\right\}+\sum_{t=1}^Tt\tr\left\{Q_{t,12}SQ_{t,12}^\top\right\}\\
   &- \tr\left\{\Lambda \left(I-\sum_{k=1}^TQ_{t,12}^\top Q_{t,12}\right)\right\}\\
    =& \inf_{Q_{t,12},t=1,\ldots,T}\;  \sum_{t=1}^T  \tr\left\{2Q_{t,12}\tilde{W}+Q_{t,12}(tS+\Lambda)Q_{t,12}^\top \right\}-\tr\{\Lambda\}
\end{align}
For this expression to be larger than $-\infty$, $S+\Lambda\succ0$ is required as otherwise we can take $Q_{t,12}=-nn^\top$ where $n$ is any vector in the null-space of $S$ and obtain that the above expression is over-bounded by
\begin{align}
    -2Tn^\top\tilde{W}n-\tr\{\Lambda\}
\end{align}
which can be made arbitrarily small by letting $\|n\|\rightarrow\infty$. For $S+\Lambda\succ0$ we have 
\begin{align}
   &{\mathcal D}(\Lambda)= \\
    =& \inf_{Q_{t,12},t=1,\ldots,T}\;  \sum_{t=1}^T  \tr\left\{\left(Q_{t,12}+\tilde{W}(tS+\Lambda)^{-1}\right)(tS+\Lambda)\left(Q_{t,12}+\tilde{W}(tS+\Lambda)^{-1}\right)\right\}\\
    &-\sum_{t=1}^T\tr\left\{\tilde{W}(tS+\Lambda)^{-1}\tilde{W}\right\}-\tr\left\{\Lambda\right\}\\
    \label{dd5}
=&-\sum_{t=1}^T\tr\left\{\tilde{W}(tS+\Lambda)^{-1}\tilde{W}\right\}-\tr\left\{\Lambda\right\}
\end{align}
where the minimum is obtained by taking 
\begin{align}
    Q_{t,12}=-\tilde{W}(tS+\Lambda)^{-1}
\end{align}
Since \eqref{d5} is strictly feasible, Slater's condition gives that the maximum of \eqref{dd5} on $\Lambda\succeq 0$ is the minimum of \eqref{d5}. We are thus led to consider the minimization problem
\begin{align}
    \inf_\Lambda\quad & \sum_{t=1}^T\tr\left\{\tilde{W}(tS+\Lambda)^{-1}\tilde{W}\right\}+\tr\left\{\Lambda\right\}\\
    \Lambda+S&\succ0\\
    \Lambda &\succeq 0
\end{align}
which by Schur complement is equivalent to 
\begin{align}
\label{dd6}
    \inf_{\Lambda,Z_t,t=1,\ldots,T}\quad \sum_{t=1}^T\tr\{Z_t\}+\tr\{\Lambda\}&\\
\text{s.t.} \;  \begin{bmatrix} Z_t &\tilde{W}\\ \tilde{W} & \Lambda+t S\end{bmatrix}&\succeq 0,\;t=1,\ldots,T\\
S+\Lambda&\succ0 \\
\Lambda&\succeq 0
\end{align}
As problem \eqref{d5} is strictly feasible, by Slater's condition the minimum of \eqref{dd6} equals the minimum of \eqref{d5}, but with opposite sign. Furthermore, as already noted, also by Slater's condition the minimum of \eqref{d5} equals the minimum of \eqref{p5}, but with opposite sign. Thus the minimum of \eqref{dd6} equals the minimum of \eqref{p5}. Let $\Lambda^*$ be the solution to \eqref{dd6}. Then in view of that \eqref{p5} is equivalent to \eqref{p6}, we see that taking $X_1=\Lambda^*$ and $X_t=0$, $t=2,\ldots,T$, is a feasible point for \eqref{p5} for which the objective function takes the same value as the minimum of \eqref{dd6}. In summary, we have shown that when \eqref{lazycond} does not hold, the solution to Problem \eqref{p5} satisfies $X_t=0$, $t=,\ldots,T$, which proves i. in the theorem. Now the condition that $X_1+S\succ0$ gives also the only if part in ii. 

It remains to obtain the bounds on the minimum objective function in iii. and iv. First we notice that \eqref{sol5} corresponds to a feasible point also when $S$ is singular, with objective function \eqref{Jmin5}. This is thus an upper bound of the minimum objective function, giving the upper bound in \eqref{SsingularJbounds}. Continuing with the case where $S$ is singular, let $\{n_k\}$ be as in iv. of the theorem and take 
\begin{align}
    Q_{t,12}=-\frac{1}{\sqrt{T}}\sum_{k=1}^{m} n_kn_k^\top
\end{align}
then 
\begin{align}
    I-\sum_{t=1}^TQ_{t,12}Q_{t,12}^\top=I-\sum_{k=1}^mn_kn_k^\top\succeq 0
\end{align}
implying that this corresponds to a feasible point of \eqref{d5}. The corresponding objective function is 
\begin{align}
    -2\sqrt{T}\sum_{k=1}^mn_k^\top\tilde{W}n_k
\end{align}
which, using that \eqref{d5} is the negative of the dual of \eqref{p5}, gives the lower bound in \eqref{SsingularJbounds}. 

Now consider that $S\succ0$ but that \eqref{lazycond} does not hold.
Take
\begin{align}
    Q_{t,12}=-\frac{c_T}{t}\tilde{W}S^{-1}
\end{align}
Then 
\begin{align}
    \sum_{t=1}^T Q_{t,12}^\top Q_{t,12}=c_T^2\beta_TS^{-1}WS^{-1}=\frac{S^{-1}WS^{-1}}{\lambda_{\max}\{S^{-1}WS^{-1}}\leq I
\end{align}
so that this corresponds to a feasible point of \eqref{d5}. The corresponding objective function is
\begin{align}
    -2c_T\alpha_T J_T^*+c_T^2\alpha_TJ_T^*=(c_T-2)c_T\alpha_TJ_T^*
\end{align}
Again using duality theory, this gives that 
\begin{align}
    (2-c_T)c_T\alpha_TJ_T^*\leq J(\{X_t^*\}_{t=1}^T)
\end{align}
which is the lower bound in iii. 

The upper bound in iii. is obtained by taking $X_t=0$, $t=2,\ldots,T$ and 
\begin{align}
    X_1=\sqrt{T}\sum_{k=1}^mn_k^\top\tilde{W}n_k\; n_kn_k^\top
    \end{align}
    and using that
    \begin{align}
        (X_1+tS)^{-1}=\frac{1}{\sqrt{T}}\sum_{k=1}^m\frac{n_kn_k^\top}{n_k^\top\tilde{W}n_k}+\frac{1}{t}\sum_k \frac{e_ke_k^\top}{\lambda_k}
    \end{align}
   holds due to the construction of $X_1$. 
This concludes the proof. \hfill\rule{2mm}{2mm}

\end{document}